\documentclass[conference,a4paper]{IEEEtran}

\usepackage[colorinlistoftodos]{todonotes}

\IEEEoverridecommandlockouts
\usepackage{cite}
\usepackage{amsmath,amssymb,amsfonts}
\usepackage{graphicx}
\usepackage{textcomp}
\usepackage{xcolor}
\usepackage{multirow}
\usepackage{array}
\usepackage{algorithm}
\usepackage{algpseudocode}
\usepackage{caption}
\def\BibTeX{{\rm B\kern-.05em{\sc i\kern-.025em b}\kern-.08em
    T\kern-.1667em\lower.7ex\hbox{E}\kern-.125emX}}
\begin{document}

\newcommand{\argmin}{\mathop{\rm argmin}}

\allowdisplaybreaks

\title{Vulnerability-Based Optimal Grid Defense Strategies for Enhancing Cyber-Physical Energy System Resilience
\thanks{This work is funded by the Deutsche Forschungsgemeinschaft (DFG, German Research Foundation), project number 360352892, priority programme DFG SPP 1984.}
}

\author{
\IEEEauthorblockN{ Eric Tönges}
\IEEEauthorblockA{\textit{Sustainable Electrical Energy Systems} \\
\textit{University of Kassel}\\
Kassel, Germany \\
eric.toenges@uni-kassel.de}
\and
\and
\IEEEauthorblockN{Martin Braun}
\IEEEauthorblockA{\textit{Sustainable Electrical Energy Systems} \\
\textit{University of Kassel, Fraunhofer IEE}\\
Kassel, Germany \\
martin.braun@iee.fraunhofer.de}
\and
\IEEEauthorblockN{Philipp Härtel}
\IEEEauthorblockA{\textit{Sustainable Electrical Energy Systems} \\
\textit{University of Kassel, Fraunhofer IEE}\\
Kassel, Germany \\
philipp.haertel@iee.fraunhofer.de}
}

\IEEEoverridecommandlockouts
\IEEEpubid{\makebox[\columnwidth]{979-8-3315-6520-6/25/\$31.00 ©2025 IEEE \hfill} \hspace{\columnsep}\makebox[\columnwidth]{ }}
\maketitle
\IEEEpubidadjcol

\begin{abstract}

An approach is proposed to identify optimal asset protection strategies based on vulnerability assessment outcomes. Traditional bilevel attacker-defender models emphasize worst-case scenarios but offer limited defensive guidance. In contrast, trilevel models introduce high computational complexity and rely on fixed network configurations. The proposed critical-components method leverages vulnerability assessment results to determine protection strategies, effectively outsourcing the upper-level defense decision. This enables adaptability to diverse network topologies, assessment techniques, and cyber-physical energy systems without the overhead of multi-level optimization. Case studies demonstrate the potential for improved system resilience across varying operational conditions.
\end{abstract}


\section*{Nomenclature}
\subsection{Sets and Indices}
\begin{tabular}{@{} l p{0.7\linewidth} @{}}
    $A_w$ & Identified critical attack scenario (CAS) with $|A_w|$ attacked components: \\
    & $A_w=\{a_{(w,1)},\ldots,a_{(w,Z^\text{max})}\}$, $a \in K\cup G$ \\
    $g \in G$ & Index of ICT-controlled generators \\
    $i \in N$ & Index of buses \\
    $k \in K$ & Index of branches \\
    $w \in W$ & Index of identified CASs
\end{tabular}
\subsection{Parameters}
\begin{tabular}{@{} l p{0.75\linewidth} @{}}
    $B_k$ & Total line susceptance of line $k$ ($\frac{\text{A}}{\text{V}}$) \\
    $P_i^{\text{dem}}$ & Active-power demand at bus $i$ (MW) \\
    $P_g^{\text{gen,max}}$ & Max. generation at generator $g$ (MW) \\
    $P_k^{\text{br,max}}$ & Max. power flow on branch $k$ (MW) \\
    $B(g)$ & Bus of generator $g$ \\
    $R(k)$ & Receiving bus of line $k$ \\
    $S(k)$ & Sending bus of line $k$ \\
    $X^{\text{max}}$ & Max. number of protected components \\
    $Z^{\text{max}}$ & Max. number of attacked components
\end{tabular}
\subsection{Decision variables}
\begin{tabular}{@{} l p{0.75\linewidth} @{}}
    $b_w$ & Equals $1$ if CAS $A_w$ is excluded, $0$ otherwise \\
    $p_i^{\text{ls}}$ & Load shedding at bus $i$ (MW) \\
    $p_g^{\text{gen}}$ & Power injection of generator $g$ (MW) \\
    $p_k^{\text{br}}$ & Power flow on line $k$ from $S(k)$ to $R(k)$ (MW) \\
    $x_k^\text{ps},x_g^\text{ict}$ & Equals $1$ if branch $k$ / generator $g$ is protected, and $0$ otherwise \\
    $y_w$ & Equals $1$ if all CASs $\{A_1,\ldots,A_w\}$ are excluded, $0$ otherwise \\
    $z_k^\text{ps},z_g^\text{ict}$ & Equals $0$ if branch $k$ / generator $g$ is attacked, $1$ otherwise \\
    $\theta_i$ & Voltage angle at bus $i$ (rad)
\end{tabular}

\section{Introduction} \label{sec:introduction}
The ongoing transformation towards net-zero emissions, decentralization, and digitalization makes distribution systems integral to power generation and supply. 
Integrating information and communication technologies (ICT) is evolving these systems into cyber-physical energy systems (CPES). 
While this facilitates coordination of new flexibility, it also introduces vulnerabilities related to high-impact-low-probability (HILP) events~\cite{Braun.2024}, triggered by deliberate attacks, such as the 2015 Ukraine blackout \cite{Liang.2017}, or extreme weather.

Enhancing the resilience of future CPES requires understanding these vulnerabilities and developing effective remedial measures. 
Mathematical bilevel optimization is a recognized approach for identifying HILP events in power systems \cite{Salmeron.2004, Zhao.2013, Abedi.2021, Ding.2017} and has been extended to CPES \cite{Castillo.2019}. 
The upper level maximizes damage to the power system under constrained attack resources, targeting components, and anticipating the grid operator’s remedial actions. 
The optimal attack represents the worst-case scenario, while critical attack scenarios (CASs) encompass the second to $n$-th best solutions. 
The lower level minimizes damage based on upper-level decisions, typically through an optimal power flow (OPF) formulation.

The literature describes mathematical trilevel models with an additional upper-level defender, allocating defense resources to protect a set of critical components from attacks \cite{Alguacil.2014, Yuan.2014, Wu.2017, Li.2024b}. 
This adds another layer of complexity to an intricate optimization problem. 
Existing approaches, such as implicit enumeration algorithms \cite{Alguacil.2014} and column-and-constraint generation algorithms (CCGA) \cite{Yuan.2014, Wu.2017, Li.2024b}, require solving the trilevel model multiple times for different defense budgets and load cases.

To our knowledge, no current approach identifies optimal protected sets based on critical solutions obtained from vulnerability assessments across varying configurations, an essential step for informing resilience-enhancing protective measures. 
We address this gap by first applying an algorithm to identify multiple CASs using bilevel optimization.
Then, we introduce a method to identify optimal protected sets that exclude CASs, minimizing load shedding in the worst remaining scenario. 
The approach extends the CCGA framework to handle a general list of CASs.
It is readily transferable to vulnerability assessments for CPES.
Our main contributions are:
\begin{itemize}
\item An integer programming formulation for identifying optimal protected components based on known vulnerabilities, with applicability to the cyber-physical domain,
\item A modular method, combinable with any vulnerability assessment approach, to determine protected-component sets across varying grid configurations and scenarios, enhancing resilience to HILP events,
\item Case studies on benchmark grids demonstrating the applicability and efficiency of the proposed approach.
\end{itemize}

The remainder is structured as follows: Section \ref{sec:bilevel} reviews existing methods for identifying CASs and introduces a mathematical formulation of bilevel network interdiction models. Section \ref{sec:trilevel} outlines the proposed method for determining optimal protected sets. Case studies on various grids are presented in Section \ref{sec:case_studies}, with results discussed in Section \ref{sec:discussion}, and the conclusion in Section \ref{sec:conclusion}.


\section{Identification of Critical Attack Scenarios} \label{sec:bilevel}
Identifying CASs in power systems and CPES is closely related to determining worst-case attacks. Various methodologies exist for detecting these solutions, often focusing on power systems, yielding lists or sets of critical components.

\subsection{Existing Approaches} \label{sec:literature_bilevel}
Identifying vulnerabilities is essential for developing remedial strategies. 
A variety of approaches exist for identifying critical components or attack vectors \cite{Abedi.2019}.
They include those based on graph theory \cite{Biswas.2021}, fault chain theory \cite{Wang.2011}, and heuristic methods focused on line loading \cite{Bier.2007}. 
Adversarial-learning-based approaches such as \cite{Fischer.2018} have also emerged in recent years.
While most existing approaches are suitable for identifying many vulnerabilities in power systems, it is crucial for resilience assessment to detect all existing vulnerabilities without overlooking any. 
Mathematical bilevel optimization models are employed to identify the worst-case attack scenario given a specific attack budget.
Building on \cite{Salmeron.2004}, various solution strategies and model extensions have emerged in recent years, for example, including line switching \cite{Zhao.2013}, cyber-physical properties \cite{Castillo.2019} and AC OPF application \cite{Abedi.2021}. 
These mathematical optimization models can identify the worst-case attack vector while incorporating the physical power flow terms.

\subsection{Bilevel Network Interdiction Models} \label{sec:model_bilevel}

Bilevel network interdiction models, particularly attacker-defender models, are crucial for identifying vulnerabilities in power systems and CPES, as they transparently incorporate physical properties of the system.
We add basic CPES features to a widely-used bilevel model that utilizes a DC OPF formulation in the lower level, as formulated in \eqref{eq:att_obj}--\eqref{eq:branch_limit}.
\begin{align}
& \max_{z} \sum_{i \in N} p_i^{\text{ls}^*} \label{eq:att_obj}\\
& \text{s.t.} \nonumber\\
& \sum_{k \in K} (1 - z_k^\text{ps}) + \sum_{g \in G} (1 - z_g^\text{ict}) \leq Z^{\text{max}} \label{eq:att_budget}\\
& z_g^\text{ict} \in \{0, 1\} \quad \forall g \in G, \quad z_k^\text{ps} \in \{0,1\} \quad \forall k \in K \label{eq:bilevel_binary}\\
& p_i^{\text{ls}^*} \in \argmin_{p^{\text{ls}}, p^{\text{gen}}, p^{\text{br}}, \theta} \left \{ \sum_{i \in N} p_i^{\text{ls}} \right \} \quad \forall i \in N \label{eq:ll_obj}\\
& \text{s.t.} \nonumber \\
& \sum_{g|B(g)=i} p_g^{\text{gen}} + p_i^{\text{ls}} - P_i^{\text{dem}} - \sum_{S(k) = i} p_k^{\text{br}} + \sum_{R(k) = i} p_k^{\text{br}} = 0 \;\nonumber \\
& \hspace{6cm} \forall i \in N \label{eq:bus_balance}\\
& p_k^{\text{br}} = z_k^\text{ps} B_k (\theta_{S(k)} - \theta_{R(k)}) \quad \forall k \in K \label{eq:branch_flow}\\
& 0 \leq p_g^{\text{gen}} \leq z_g^\text{ict} P_g^{\text{gen,max}} \quad \forall g \in G \label{eq:gen_lim}\\
& 0 \leq p_i^{\text{ls}} \leq P_i^{\text{dem}} \quad \forall i \in N  \label{eq:dem_lim}\\
& -P_k^{\text{br,max}} \leq p_k^{\text{br}} \leq P_k^{\text{br,max}} \quad \forall k \in K \label{eq:branch_limit}
\end{align}
The upper level in \eqref{eq:att_obj} maximizes lost load by disabling components through binary attack variables defined in \eqref{eq:bilevel_binary}, adhering to the attack budget, i.e. the maximum number of attacked components, in \eqref{eq:att_budget}. 
Attacks can be physical attacks targeting branches (indicated by $z_k^\text{ps}$) and cyber-attacks disconnecting ICT-controlled generators ($z_g^\text{ict}$).
The lower level minimizes lost load in \eqref{eq:ll_obj} based on upper-level decisions. 
The bus power balance and active power flow on branches are represented in \eqref{eq:bus_balance} and \eqref{eq:branch_flow}, respectively. 
Inequalities \eqref{eq:gen_lim}--\eqref{eq:branch_limit} limit generation, load shedding, and maximum branch flow. 

\subsection{Identification Beyond the Worst Case} \label{sec:extension_bilevel}

To identify CASs rather than only one worst-case attack scenario, constraint \eqref{eq:bilevel_add} is added for each identified attack scenario $A_w=\{a_{(w,1)},…,a_{(w,Z^\text{max})}\}$ with $a\in K\cup G$:
\begin{align}
& \sum_{k \in A_w} z_{k}^\text{ph} + \sum_{g \in A_w} z_{g}^\text{ict} \geq 1 \quad \forall w \in W \,.\label{eq:bilevel_add}
\end{align}
This ensures that in each previously identified attack scenario, at least one physically or cyber-attacked component remains unattacked, excluding all prior solutions and enabling the detection of an additional CAS to set $W$, similar to the methodology reported in \cite{Ding.2017}.

\section{Identification of Asset Protection Strategies} \label{sec:trilevel}
Identifying vulnerabilities and CASs is a critical first step in enhancing resilience in power systems and CPES. Equally important is determining how to protect the system against these threats. Prior work \cite{Yuan.2014, Wu.2017, Alguacil.2014, Li.2024b} addresses the optimal allocation of protective resources to minimize worst-case impacts.
Trilevel defender-attacker-defender models extend the bilevel attacker-defender models from Section \ref{sec:model_bilevel} by identifying the best set of components to protect.


\subsection{Trilevel Defender-Attacker-Defender Models} \label{sec:existing_trilevel}

The trilevel defender-attacker-defender model is obtained by adding two modifications to the bilevel model in \eqref{eq:att_obj}--\eqref{eq:branch_limit}. First, an additional upper level is introduced in \eqref{eq:trilevel_obj}--\eqref{eq:trilevel_binary}:
\begin{align}
& \min_x \sum_{i \in N} p_i^{\text{ls}^*} \label{eq:trilevel_obj}\\
& \text{s.t.} \nonumber \\
& \sum_{k \in K} x_k^\text{ps} + \sum_{g \in G} x_g^\text{ict} \leq X^{\text{max}} \\
& x_g^\text{ict} \in \{0, 1\} \quad \forall g \in G, \quad x_k^\text{ps} \in \{0, 1\} \quad \forall k \in K \label{eq:trilevel_binary}
\end{align}
The defender's objective is to minimize lost load by allocating the maximum defense budget, i.e. the maximum number of components protected from attacks, to attackable components. The binary variable $x_k^\text{ps}$ ($x_g^\text{ict}$) equals $1$ if branch $k$ (generator $g$) is excluded from attacks, e.g. through physical grid hardening or cyber-security measures, and $0$ otherwise. The second modification adds \eqref{eq:trilevel_bilevel} to the attacker model in \eqref{eq:att_obj}--\eqref{eq:bilevel_binary}:
\begin{align}
& \quad z_g^\text{ict} \geq x_g^\text{ict} \quad \forall g \in G,\quad z_k^\text{ps} \geq x_k^\text{ps} \quad \forall k \in K \,. \label{eq:trilevel_bilevel}
\end{align}
This ensures protected components cannot be taken out of service.
The trilevel structure adds complexity to the existing bilevel attacker-defender model.
Various approaches for solving these models for power systems have been proposed, including a branching-based enumeration method \cite{Alguacil.2014} and a CCGA that decomposes the trilevel model into a master problem and a subproblem, similar to Bender’s decomposition \cite{Yuan.2014, Wu.2017}. 
The CCGA has been shown to outperform the enumeration method \cite{Yuan.2014} and is extended to CPES in \cite{Li.2024b}.

However, integrated trilevel models require predefined configurations for $X^\text{max}$ and $Z^\text{max}$, and input parameters related to load and generation, limiting identified solutions to  specified grid configurations. 
Changes in load cases, attack budgets, or protection budgets necessitate re-running the trilevel model. 
Moreover, most existing approaches utilize a DC OPF in the lower level, incorporating an AC OPF adds further complexity to the overall model.

\subsection{Vulnerability-Based Approaches} \label{sec:new_approaches}

\begin{figure}[tb]
    \centering
    \includegraphics[width=\linewidth]{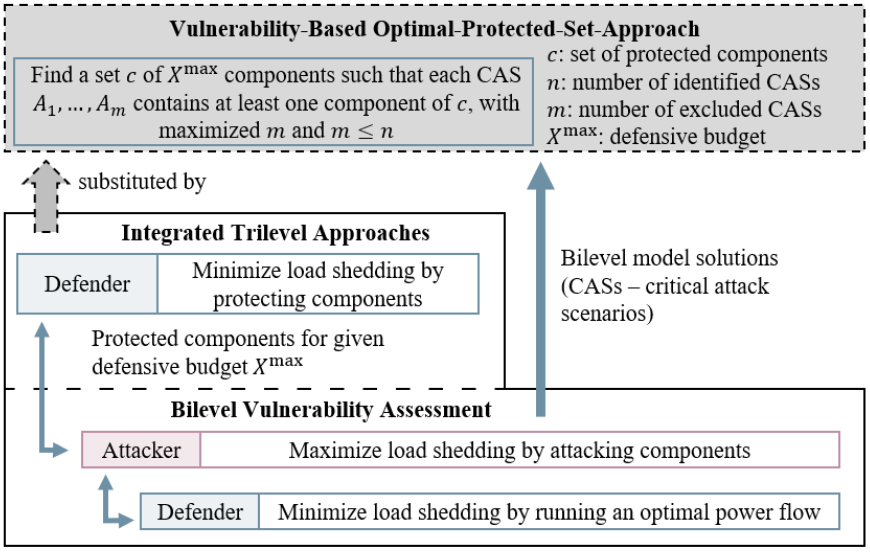}
    \caption{Schematic comparison of integrated trilevel defender-attacker-defender models and the proposed vulnerability-based approach to identify optimal asset protection strategies.}
    \label{fig:approach_overview}
\end{figure}

In order to overcome the drawbacks of integrated trilevel optimization models, we propose an approach for identifying the optimal set of protected components based on a precomputed list of CASs, as illustrated in Fig. \ref{fig:approach_overview}. 
While the pre-identification may affect runtime, it provides valuable insights into grid-related vulnerabilities and accommodates various vulnerability identification methods beyond mathematical bilevel models, such as mentioned in Section \ref{sec:literature_bilevel}.
For determining the optimal set of protected components, an integer optimization model, similar to established CCGA approaches is given in \eqref{eq:intopt_obj}--\eqref{eq:intopt_last}.

\begin{align}
& \max_{b,x,y} \sum_{w \in W} y_w \label{eq:intopt_obj} \\ 
& \text{s.t.} \nonumber \\
& \sum_{k \in K} x_k^\text{ps} + \sum_{g \in G} x_g^\text{ict} \leq X^{\text{max}} \label{eq:intopt_budget}\\
& x_g^\text{ict} \in \{0,1\} \quad \forall g \in G, \quad x_k^\text{ps} \in \{0,1\} \quad \forall k \in K \label{eq:intopt_binary1}\\
& y_w, b_w \in \{0, 1\} \quad \forall w \in W \label{eq:intopt_binary2} \\
& b_w \leq \sum_{k \in A_w} x_k^\text{ps} + \sum_{g \in A_w} x_g^\text{ict} \leq |A_w| b_w \quad \forall w \in W \label{eq:intopt_limits2} \\
& y_w \leq b_w \quad \forall w \in W \label{eq:intopt_consec1}\\
& y_w \geq y_{w+1} \quad \forall w = \{1, \ldots, |W|-1\} \label{eq:intopt_last}
\end{align}

The objective in \eqref{eq:intopt_obj} maximizes the number of consecutively excluded CASs, indicated by the binary variable $y_w$, starting from the worst-case scenario.
The binary variables defined in \eqref{eq:intopt_binary1}--\eqref{eq:intopt_binary2} are used to represent protection decisions: $x_k^\text{ps}$ for branches protected from physical attacks, $x_g^\text{ict}$ for generators protected from cyber-attacks, $b_w$ and $y_w$ indicating exclusion of scenario $w$ and of scenarios \{$1, ..., w$\}, respectively. 
The protection budget constraint is defined in \eqref{eq:intopt_budget}.
Constraint \eqref{eq:intopt_limits2} ensures $b_w=1$ if at least one component in the attack vector $A_w$ is protected, while \eqref{eq:intopt_consec1}--\eqref{eq:intopt_last} guarantee that $y_w$ equals $1$ only if the considered attack and all previous attacks are excluded.

Given a set of CASs with associated lost load, the approach identifies the optimal solution candidate set.
If multiple candidate sets yield an optimal solution, the optimization may return only one, but additional next-best solutions can be provided depending on the applied solver and solver options.

The proposed model operates independently of the vulnerability assessment method employed during the identification stage. Consequently, no adaptations or modifications to the model are necessary, regardless of whether heuristics, learning-based methods~\cite{Lai.2023}, or mathematical approaches that focus on AC OPF or more detailed ICT representations are utilized, provided that a list of CASs is given.

\section{Case Studies} \label{sec:case_studies}
\subsection{Case Study Design} \label{sec:study_setup}

The applicability of the proposed methodology is demonstrated through case studies on benchmark grids of varying sizes, drawn from pandapower~\cite{Thurner.2018b} and SimBench~\cite{Meinecke.2020} as shown in Table~\ref{tab:grids_overview}.
\begin{table}[tb]
\caption{Overview of the considered grids and their configurations for the case studies.}
\label{tab:grids_overview}
\resizebox{\linewidth}{!}{%
\begin{tabular}{|c|c|c|c|c|}
\hline
Grid description & Configurations & \#buses & \#branches & \#gens \\ \hline
IEEE 9-bus \cite{Thurner.2018b} & standard & 9 & 9 & 2\\ \hline
\multirow{2}{*}{CIGRE MV (with DER) \cite{Thurner.2018b}} & open switches & \multirow{2}{*}{15} & \multirow{2}{*}{17} & \multirow{2}{*}{13} \\
 & closed switches &  &  & \\ \hline
IEEE 30-bus \cite{Thurner.2018b} & standard & 30 & 41 & 5 \\ \hline
\multirow{2}{*}{\begin{tabular}[c]{@{}c@{}}SimBench\\ 1-HV-urban--0-no\_sw \cite{Meinecke.2020}\end{tabular}} & high-load case & \multirow{2}{*}{82} & \multirow{2}{*}{116} & \multirow{2}{*}{98} \\
& low-load case &  &  & \\ \hline
\multirow{2}{*}{\begin{tabular}[c]{@{}c@{}}SimBench\\ 1-HV-urban--2-no\_sw \cite{Meinecke.2020}\end{tabular}} & high-load case & \multirow{2}{*}{120} & \multirow{2}{*}{154} & \multirow{2}{*}{118} \\ & low-load case &  &  & \\ \hline
\end{tabular}%
}
\end{table}
The CIGRE MV grid is analyzed both as a radial structure with open switches and as a version with closed switches. 
Each SimBench grid is evaluated in two time steps: one with high load and low generation (high-load case), and another with low load and high generation (low-load case).
In the showcased SimBench grids, which involve a high number of distributed and controllable generators, low-load cases with significant feed-in are anticipated to exhibit less criticality.
However, depending on the grid configuration and the vulnerability assessment approach employed, these scenarios may also reveal critical vulnerabilities.
Furthermore, considering two distinct load cases in this study demonstrates the applicability of the proposed defense approach across a range of grid configurations and load conditions.
Several modeling assumptions are applied consistently across all grids: storage units and shunt elements are excluded, transformer taps are fixed at neutral positions, loads and generators are treated as continuously controllable, and minimum generator output is set to zero to avoid dualization infeasibilities \cite{Bienstock.2010}.

\begin{figure}[tb]
    \centering
    \includegraphics[width=\linewidth]{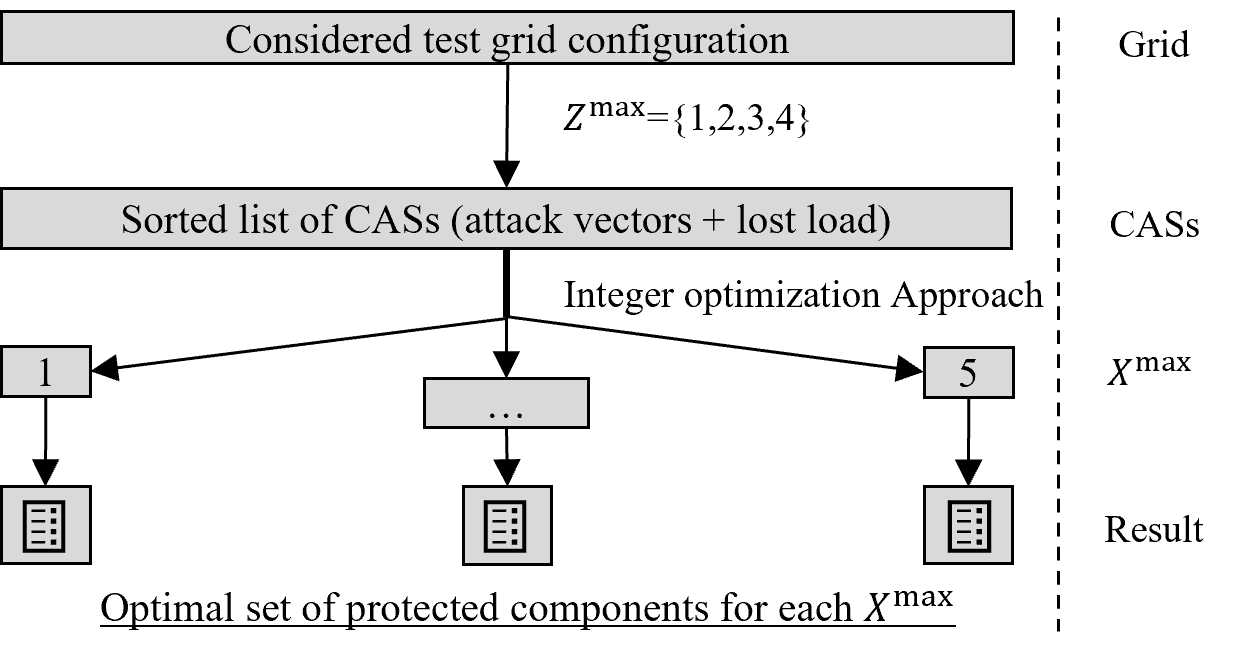}
    \caption{Simulation procedure applied in the case studies.}
    \label{fig:simulation_overview}
\end{figure}

For each grid configuration, the bilevel model described in Section \ref{sec:model_bilevel}, along with the extension from Section \ref{sec:extension_bilevel}, is applied to identify CASs and the associated lost load. 
The simulation procedure is illustrated in Fig.~\ref{fig:simulation_overview}.
A total of 216 to 526 CASs are generated per configuration, considering up to $Z^{\text{max}}=4$ attackable components. 
The integer optimization approach is then executed with a protection budget of $X^{\text{max}}=\{1,2,3,4,5\}$. 
Although these settings are assumed in the case studies, $X^{\text{max}}$ and $Z^{\text{max}}$ can be selected freely, depending on the goals of the analysis.
For each case, the optimal protected set and corresponding runtime are recorded.
The calculation of protected sets is repeated for a combined list of CASs across both time steps in each SimBench grid and for both CIGRE configurations, demonstrating the adaptability of the approach.

\subsection{Case Study Results} \label{sec:study_results}

\begin{figure}[tb]
    \centering
    \includegraphics[width=\linewidth]{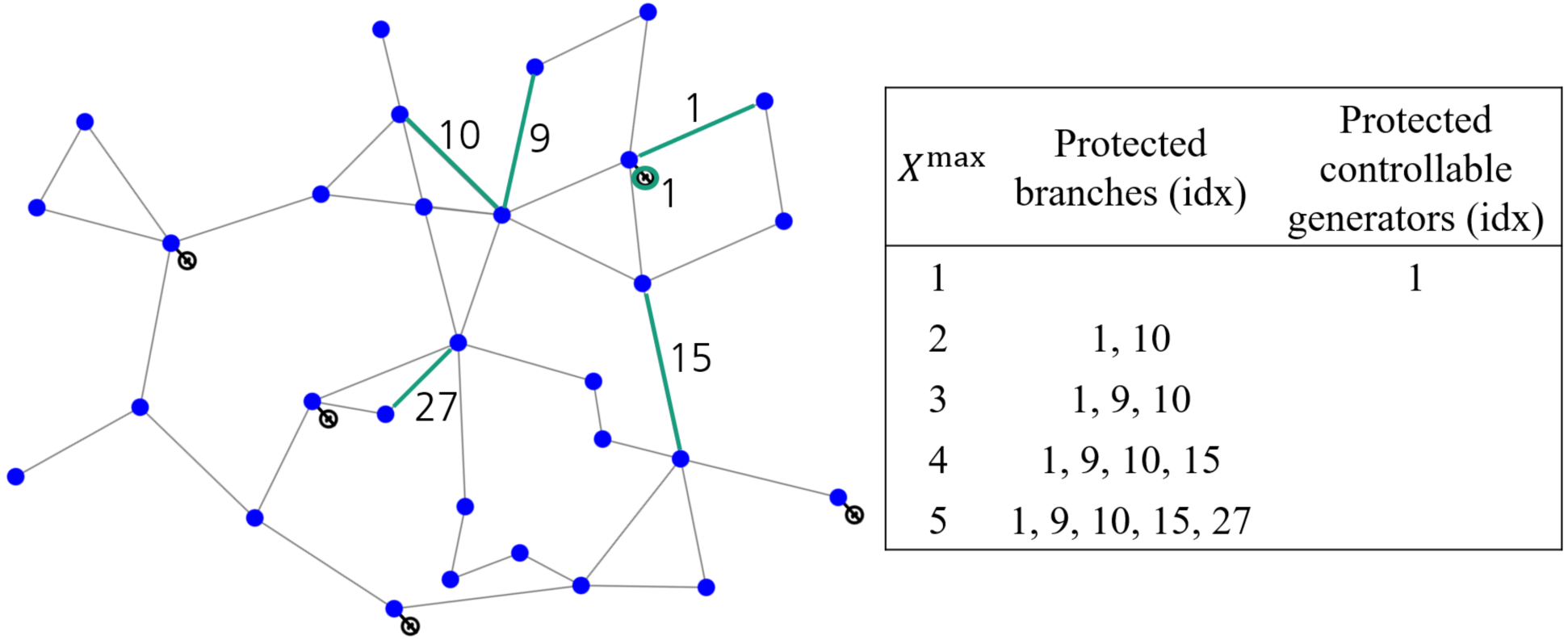}
    \caption{Optimal protected sets of components with varying $X^\text{max}$ in the IEEE 30-bus grid.}
    \label{fig:ieee30_results}
\end{figure}

Fig. \ref{fig:ieee30_results} illustrates the topology of the IEEE 30-bus system, highlighting the protected components for various protection budgets $X^\text{max}=\{1,2,3,4,5\}$, obtained with the integer optimization model.
A key feature of defender-attacker-defender models visible in Fig.~\ref{fig:ieee30_results} is that components included in the optimal set for $X^\text{max}=1$ may differ from those for higher budgets.
Specifically, the optimal protection decision for $X^\text{max}=1$ selects generator $1$, which is not part of the optimal set for $X^\text{max}=2$, where branches $1$ and $10$ are chosen instead.
A similar pattern is observed in the bilevel models with varying attacker budgets $Z^\text{max}$, where a CAS of size two may contain different components than a CAS of size three.
These results expose the benefit of mathematical optimization methods, such as multilevel optimization, for systematically identifying CASs and protection sets under varying budget constraints.

\begin{figure}[tb]
    \centering
    \includegraphics[width=\linewidth]{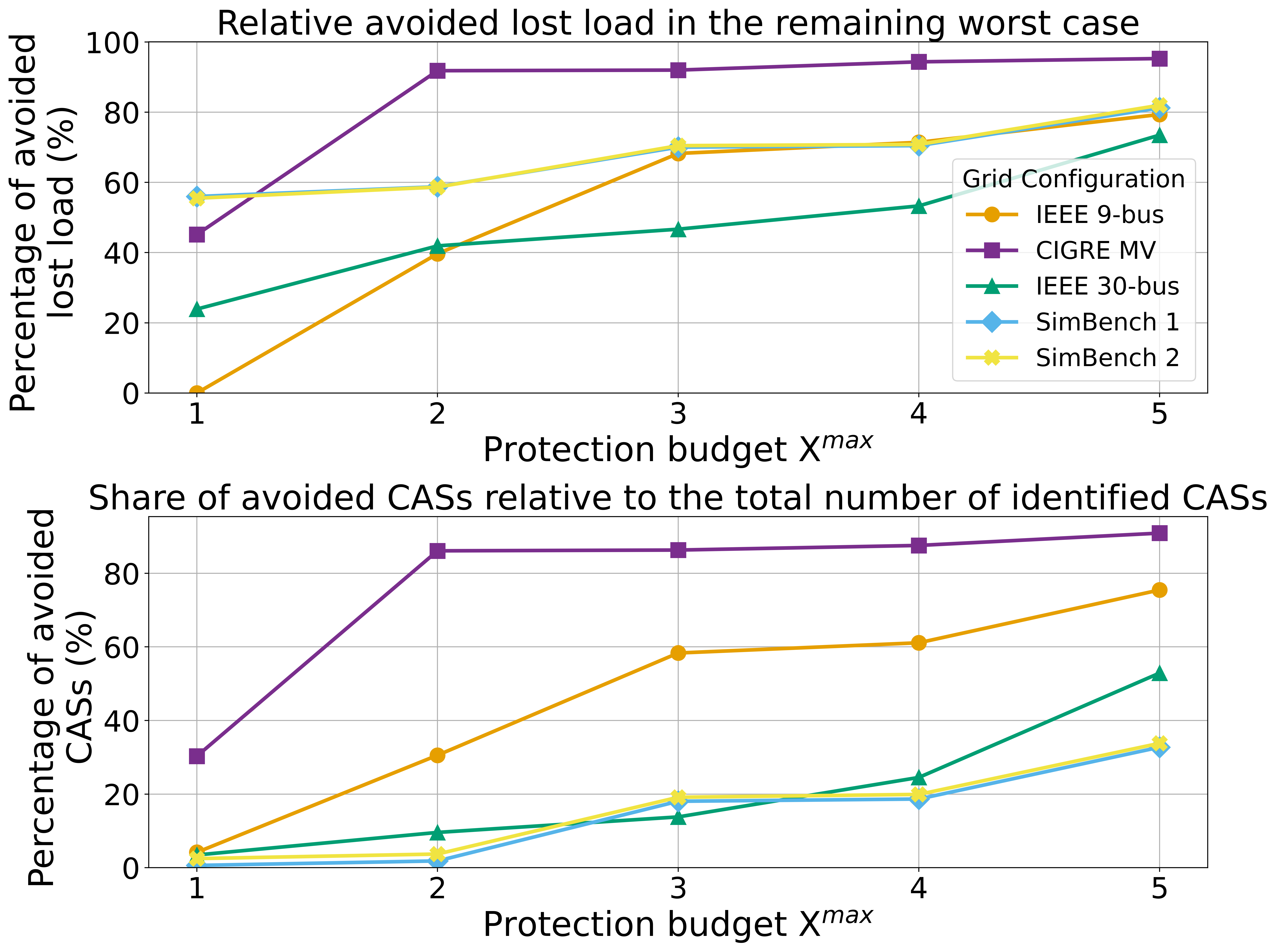}
    \caption{Effect of varying $X^\text{max}$ on the reduction of remaining worst-case lost load and the total number of excluded CASs.}
    \label{fig:protection_results}
\end{figure}

The impact of varying $X^\text{max}$ on the excluded CASs and the avoided worst-case lost load for each grid is shown in Fig.~\ref{fig:protection_results}. 
For the CIGRE grid and both SimBench grids, the results of the combined analysis for both grid configurations are used.
The first subplot illustrates the percentage of avoided potential lost load in the remaining worst-case scenario. 
The second subplot depicts the reduction in the absolute number of CASs remaining unprotected by the optimal set, compared to the initial number of identified CASs.
For example, the lost load in the worst unprotected scenario for the IEEE 9-bus grid with $X^\text{max}=2$ is approximately 40\,\% lower compared to the case without any protection budget. 
Around 30\,\% of the worst identified CASs for this grid are prevented by the optimal protected set.
In the case of the IEEE 30-bus system shown in Fig. \ref{fig:ieee30_results}, protecting generator $1$ excludes approximately 5\,\% of the worst CASs but reduces the remaining worst-case lost load by more than 20\,\%. 
The protection of the branch set \{1, 10\} reduces the remaining worst-case lost load by around 40\,\% compared to the scenario without protection, while excluding around 10\,\% of the worst CASs.

For grids connected to an external grid (CIGRE and SimBench), transformers are the most critical components regarding lost load. 
A budget of $X^\text{max}=1$ significantly impacts the avoided lost load, reducing it by approximately half for each configuration. 
In contrast, meshed IEEE grids require a higher protection budget to achieve comparable effects. 
Notably, the share of avoided lost load increases more than the share of avoided CASs, indicating that even a small number of identified CASs can lead to substantial damage reduction, while the contribution of increasing $X^\text{max}$ is also evident. 
All calculations were performed using a laptop with an Intel Core i7, 2.8\,GHz processor, 16\,GB RAM, and Gurobi 12 \cite{GurobiOptimization.2025}, achieving runtimes around one second or lower.
Scalability across the considered test grids, ranging from 9 to 120 buses, is demonstrated, with runtimes expected to remain manageable for larger, realistic energy systems.
The main scalability constraint arises in the preliminary step of vulnerability identification.
Thanks to the computational efficiency of our method, scalability challenges are limited to this stage and are not compounded by the defense optimization process itself.
Moreover, the approach is adaptable to alternative, more computationally efficient vulnerability assessment methods beyond mathematical bilevel optimization, enabling its application to larger systems for determining optimal defense strategies.
Overall, the case studies demonstrate the applicability and scalability of the approach across various grid configurations and topologies, and suitability for application in larger real-world grids is expected.

\section{Discussion} \label{sec:discussion}
Our methodology aims to enhance the applicability of optimal defensive resource allocation in trilevel defender-attacker-defender models. 
While similar to CCGA algorithms in integrated trilevel models, the approaches differ in features and use cases.
These differences are summarized in Table~\ref{tab:comparison_app}.
\begin{table}[tb]
\caption{Features and use cases of existing integrated trilevel approaches and the presented vulnerability-based approach.}
\label{tab:comparison_app}
\resizebox{\linewidth}{!}{%
\begin{tabular}{|c|c|c|}
\hline
 & \begin{tabular}[c]{@{}c@{}}Integrated Trilevel \\ Optimization\end{tabular} & \begin{tabular}[c]{@{}c@{}}Presented\\ Methodology\end{tabular} \\ \hline
Global optimality & $\checkmark$ & $\checkmark$ \\
Different topological configurations & $\times$ & $\checkmark$ \\
Different load cases & $\times$ & $\checkmark$ \\
Different attacker and defender budgets & $(\checkmark)$ & $\checkmark$ \\
\begin{tabular}[c]{@{}c@{}}Efficiency for a small number of\\ configurations\end{tabular} & $\checkmark$ & $\times$ \\
\begin{tabular}[c]{@{}c@{}}Applicability to other vulnerability\\ assessment approaches\end{tabular} & $\times$ & $\checkmark$ \\
Transparency regarding vulnerabilities & $\times$ & $\checkmark$ \\ \hline
\end{tabular}%
}
\end{table}

The advantages of our integer optimization approach include flexibility and adaptability. 
In real-world applications, considering varying load cases and topological configurations is crucial for capturing vulnerabilities across operational states, and considering different attacker and defender budgets is important for investment decisions. 
For each of these configurations, integrated approaches require multiple re-runs of the trilevel model, while our methodology operates on a list of CASs.

Although pre-identifying CASs may reduce computational efficiency when only a few configurations are analyzed, it significantly improves transparency. 
Grid operators can more easily interpret and accept defensive strategies when the vulnerabilities are explicitly documented.
Moreover, our approach also accommodates diverse vulnerability identification methods, requiring only a list of CASs, enabling the integration of learning-based methodologies, heuristics, and ICT-aware vulnerability assessments in CPES.

This flexibility is valuable even with regard to bilevel optimization settings.
While bilevel models such as \eqref{eq:att_obj}--\eqref{eq:branch_limit} are widely used in integrated trilevel approaches, a more accurate AC OPF can result in different worst-case attack scenarios and higher potential damages compared to DC OPF calculations \cite{Abedi.2021}. 
Also, integrating mixed-integer decisions in the lower level, which requires dedicated solution approaches \cite{Heid.2025}, and a more detailed representation of the ICT properties of CPES both add additional complexity to the attacker-defender models.
For these cases, our approach helps mitigate an even higher complexity associated with integrated trilevel optimization models.

\section{Conclusion} \label{sec:conclusion}
This work presented a methodology based on integer optimization for identifying optimal defense resource allocation in the trilevel defender-attacker-defender network interdiction problem, building on results from vulnerability assessments.
For broad applicability, the approach is designed to be independent of the specific optimization model used for vulnerability identification. 
It focuses on selecting component protections that minimize worst-case lost load, can be easily applied to cyber-physical energy systems, and delivers optimal results within one second for test cases up to 154 branches and 120 buses.
The method is expected to scale efficiently to real-world network sizes.

The methodology enhances resilience in power systems and cyber-physical energy systems by providing a flexible framework for mitigating high-impact, low-probability threats.
Future work includes the incorporation of more advanced vulnerability assessments, such as AC OPF-based bilevel models and data-driven techniques. 
Further research may also explore vulnerabilities and defense strategies across varying load cases and network configurations and include detailed cyber-physical properties.

\section*{Acknowledgement}
The Fraunhofer AI chatbot FhGenie \cite{Weber.2024}, based on GPT 4o Mini, was used to support linguistic refinement. All suggestions were carefully reviewed and edited by the authors to ensure that the technical content remained unchanged.

\bibliographystyle{IEEEtran}
\bibliography{literature}

\begin{thebibliography}{10}
\providecommand{\url}[1]{#1}
\csname url@samestyle\endcsname
\providecommand{\newblock}{\relax}
\providecommand{\bibinfo}[2]{#2}
\providecommand{\BIBentrySTDinterwordspacing}{\spaceskip=0pt\relax}
\providecommand{\BIBentryALTinterwordstretchfactor}{4}
\providecommand{\BIBentryALTinterwordspacing}{\spaceskip=\fontdimen2\font plus
\BIBentryALTinterwordstretchfactor\fontdimen3\font minus \fontdimen4\font\relax}
\providecommand{\BIBforeignlanguage}[2]{{%
\expandafter\ifx\csname l@#1\endcsname\relax
\typeout{** WARNING: IEEEtran.bst: No hyphenation pattern has been}%
\typeout{** loaded for the language `#1'. Using the pattern for}%
\typeout{** the default language instead.}%
\else
\language=\csname l@#1\endcsname
\fi
#2}}
\providecommand{\BIBdecl}{\relax}
\BIBdecl

\bibitem{Braun.2024}
M.~Braun, C.~Gruhl, C.~A. Hans, P.~H{\"a}rtel, C.~Scholz, B.~Sick, M.~Siefert, F.~Steinke, O.~Stursberg, and S.~{Wende-von Berg}, ``{Predictions and Decision Making for Resilient Intelligent Sustainable Energy Systems},'' in \emph{2024 IEEE PES Innovative Smart Grid Technologies Europe (ISGT EUROPE)}.\hskip 1em plus 0.5em minus 0.4em\relax IEEE, 2024, pp. 1--5.

\bibitem{Liang.2017}
G.~Liang, S.~R. Weller, J.~Zhao, F.~Luo, and Z.~Y. Dong, ``{The 2015 Ukraine Blackout: Implications for False Data Injection Attacks},'' \emph{IEEE Transactions on Power Systems}, vol.~32, no.~4, pp. 3317--3318, 2017.

\bibitem{Salmeron.2004}
J.~Salmeron, K.~Wood, and R.~Baldick, ``{Analysis of Electric Grid Security Under Terrorist Threat},'' \emph{IEEE Transactions on Power Systems}, vol.~19, no.~2, pp. 905--912, 2004.

\bibitem{Zhao.2013}
L.~Zhao and B.~Zeng, ``{Vulnerability Analysis of Power Grids With Line Switching},'' \emph{IEEE Transactions on Power Systems}, vol.~28, no.~3, pp. 2727--2736, 2013.

\bibitem{Abedi.2021}
A.~Abedi, M.~R. Hesamzadeh, and F.~Romerio, ``{An ACOPF-based bilevel optimization approach for vulnerability assessment of a power system},'' \emph{International Journal of Electrical Power {\&} Energy Systems}, vol. 125, p. 106455, 2021.

\bibitem{Ding.2017}
T.~Ding, C.~Li, C.~Yan, F.~Li, and Z.~Bie, ``{A Bilevel Optimization Model for Risk Assessment and Contingency Ranking in Transmission System Reliability Evaluation},'' \emph{IEEE Transactions on Power Systems}, vol.~32, no.~5, pp. 3803--3813, 2017.

\bibitem{Castillo.2019}
A.~Castillo, B.~Arguello, G.~Cruz, and L.~Swiler, ``{Cyber-Physical Emulation and Optimization of Worst-Case Cyber Attacks on the Power Grid},'' in \emph{2019 Resilience Week (RWS)}.\hskip 1em plus 0.5em minus 0.4em\relax IEEE, 2019, pp. 14--18.

\bibitem{Alguacil.2014}
N.~Alguacil, A.~Delgadillo, and J.~M. Arroyo, ``{A trilevel programming approach for electric grid defense planning},'' \emph{Computers {\&} Operations Research}, vol.~41, pp. 282--290, 2014.

\bibitem{Yuan.2014}
W.~Yuan, L.~Zhao, and B.~Zeng, ``{Optimal power grid protection through a defender--attacker--defender model},'' \emph{Reliability Engineering {\&} System Safety}, vol. 121, pp. 83--89, 2014.

\bibitem{Wu.2017}
X.~Wu and A.~J. Conejo, ``{An Efficient Tri-Level Optimization Model for Electric Grid Defense Planning},'' \emph{IEEE Transactions on Power Systems}, vol.~32, no.~4, pp. 2984--2994, 2017.

\bibitem{Li.2024b}
P.~{Li et al.}, ``{A Defense Planning Model for a Power System Against Coordinated Cyber-Physical Attack},'' \emph{Protection and Control of Modern Power Systems}, vol.~9, no.~5, pp. 84--95, 2024.

\bibitem{Abedi.2019}
A.~Abedi, L.~Gaudard, and F.~Romerio, ``{Review of major approaches to analyze vulnerability in power system},'' \emph{Reliability Engineering {\&} System Safety}, vol. 183, pp. 153--172, 2019.

\bibitem{Biswas.2021}
R.~S. Biswas, A.~Pal, T.~Werho, and V.~Vittal, ``{A Graph Theoretic Approach to Power System Vulnerability Identification},'' \emph{IEEE Transactions on Power Systems}, vol.~36, no.~2, pp. 923--935, 2021.

\bibitem{Wang.2011}
A.~Wang, Y.~Luo, G.~Tu, and P.~Liu, ``{Vulnerability Assessment Scheme for Power System Transmission Networks Based on the Fault Chain Theory},'' \emph{IEEE Transactions on Power Systems}, vol.~26, no.~1, pp. 442--450, 2011.

\bibitem{Bier.2007}
V.~M. Bier, E.~R. Gratz, N.~J. Haphuriwat, W.~Magua, and K.~R. Wierzbicki, ``{Methodology for identifying near-optimal interdiction strategies for a power transmission system},'' \emph{Reliability Engineering {\&} System Safety}, vol.~92, no.~9, pp. 1155--1161, 2007.

\bibitem{Fischer.2018}
L.~Fischer, J.-M. Memmen, E.~M. Veith, and M.~Tr{\"o}schel, ``{Adversarial Resilience Learning - Towards Systemic Vulnerability Analysis for Large and Complex Systems}.''

\bibitem{Lai.2023}
\BIBentryALTinterwordspacing
Z.~Lai, K.~Wei, Y.~Fu, P.~H\"{a}rtel, and F.~Heide, ``{$\nabla$-Prox: Differentiable Proximal Algorithm Modeling for Large-Scale Optimization},'' \emph{ACM Trans. Graph.}, vol.~42, no.~4, Jul. 2023. [Online]. Available: \url{https://doi.org/10.1145/3592144}
\BIBentrySTDinterwordspacing

\bibitem{Thurner.2018b}
L.~{Thurner et al.}, ``{Pandapower---An Open-Source Python Tool for Convenient Modeling, Analysis, and Optimization of Electric Power Systems},'' \emph{IEEE Transactions on Power Systems}, vol.~33, no.~6, pp. 6510--6521, 2018.

\bibitem{Meinecke.2020}
S.~{Meinecke et al.}, ``{SimBench---A Benchmark Dataset of Electric Power Systems to Compare Innovative Solutions Based on Power Flow Analysis},'' \emph{Energies}, vol.~13, no.~12, p. 3290, 2020.

\bibitem{Bienstock.2010}
D.~Bienstock and A.~Verma, ``{The $N-k$ Problem in Power Grids: New Models, Formulations, and Numerical Experiments},'' \emph{SIAM Journal on Optimization}, vol.~20, no.~5, pp. 2352--2380, 2010.

\bibitem{GurobiOptimization.2025}
\BIBentryALTinterwordspacing
L.~L. {Gurobi Optimization}, ``{Gurobi Optimizer Reference Manual},'' 2025. [Online]. Available: \url{https://www.gurobi.com}
\BIBentrySTDinterwordspacing

\bibitem{Heid.2025}
J.~Heid, N.~Bornhorst, E.~T\"{o}nges, P.~H\"{a}rtel, D.~Mende, and M.~Braun, ``{A Computationally Efficient Method for Solving Mixed-Integer AC Optimal Power Flow Problems},'' in \emph{2025 IEEE Kiel PowerTech}.\hskip 1em plus 0.5em minus 0.4em\relax IEEE, 2025.

\bibitem{Weber.2024}
I.~Weber, H.~Linka, D.~Mertens, T.~Muryshkin, H.~Opgenoorth, and S.~Langer, ``Fhgenie: a custom, confidentiality-preserving chat ai for corporate and scientific use,'' in \emph{2024 IEEE 21st International Conference on Software Architecture Companion (ICSA-C)}.\hskip 1em plus 0.5em minus 0.4em\relax IEEE, 2024, pp. 26--31.

\end{thebibliography}

\end{document}